\documentclass{amsart}
\usepackage{latexsym,amsmath,amssymb}
\theoremstyle{definition}
\theoremstyle{remark}

\numberwithin{equation}{section}
\begin{document}

\title[MULTIPLICATION AND COMPOSITION OPERATORS BETWEEN TWO DIFFERENT ORLICZ SPACES]
{MULTIPLICATION AND COMPOSITION OPERATORS BETWEEN TWO DIFFERENT
ORLICZ SPACES}

\author{\bf Y. Estaremi}
\address{ Y. Estaremi}
\address{ Department of Mathematics, Payame Noor University (PNU), P. O. Box: 19395-3697, Tehran- Iran}
 \email{estaremi@gmail.com}



\begin{abstract}
In this paper we consider composition operator $C_{\varphi}$
generated by nonsingular  measurable transformation
$\varphi:\Omega\rightarrow \Omega$ and multiplication operator
$M_{u}$ generated by measurable function $u:\Omega\rightarrow
\mathbb{C}$ between two different Orlicz spaces
$L^{\Phi_1}(\Omega,\Sigma, \mu)$ and $L^{\Phi_2}(\Omega, \Sigma,
\mu)$, then we investigate boundedness, compactness and essential
norm of multiplication and composition operators in term of
properties of the mapping $\varphi$, the function $u$ and the
measure space $(\Omega, \Sigma, \mu)$.
\end{abstract}
\maketitle

\section{\sc\bf Introduction and Preliminaries}

Let $\Phi:\mathbb{R}\rightarrow\mathbb{R}^{+}$ be a continuous
convex function such that

(1)$\Phi(x)=0$ if and only if $x=0$.

(2) $\lim_{x\rightarrow\infty}\Phi(x)=\infty$.

(3) $\lim_{x\rightarrow\infty}\frac{\Phi(x)}{x}=\infty$.

 The convex function $\Phi$ is called Young's function. With each
Young's function $\Phi$, one can associate another convex function
$\Psi:\mathbb{R}\rightarrow\mathbb{R}^{+}$ having similar
properties, which is defined by
$$\Psi(y)=\sup\{x|y|-\Phi(x):x\geq0\}, \ \ y\in\mathbb{R}.$$
The convex function $\Psi$ is called complementary Young function
to $\Phi$. A Young function $\Phi$ is said to satisfy the
$\bigtriangleup_{2}$ condition (globally) if $\Phi(2x)\leq
k\Phi(x), \ x\geq x_{0}\geq0 (x_{0}=0)$ for some constant $k>0$.

If $\Phi$ is a Young function, then the set of $\Sigma-$measurable functions
 $$L^{\Phi}(\Sigma)=\{f:\Omega\rightarrow \mathbb{C}:\exists k>0, \int_{\Omega}\Phi(k|f|)d\mu<\infty\}$$
 is a Banach space, with respect to the norm $N_{\Phi}(f)=\inf\{k>0:\int_{\Omega}\Phi(\frac{f}{k})d\mu\leq1\}$.
 $(L^{\Phi}(\Sigma), N_{\Phi}(.))$ is called Orlicz space. The
 usual convergence in the orlicz space $L^{\Phi}(\Sigma)$ can be
 introduced in term of the orlicz norm $N_{\Phi}(.)$ as
 $u_n\rightarrow u$ in $L^{\Phi}(\Sigma)$ means $N_{\Phi}(u_n-u)\rightarrow
 0$. Also, a sequence $\{u_n\}^{\infty}_{n=1}$ in
 $L^{\Phi}(\Sigma)$ is said to converges in $\Phi$-mean to $u\in
 L^{\Phi}(\Sigma)$, if

  $$\lim_{n\rightarrow\infty}
 I_{\Phi}(u_n-u)=\lim_{n\rightarrow\infty}\int_{\Omega}\Phi(|u_n-u|)d\mu=0.$$

 Let $\Omega=(\Omega, \Sigma, \mu)$ be a $\sigma$-finite
complete measure space and let $\varphi:\Omega\rightarrow \Omega$
be a measurable transformation, that is, $\varphi^{-1}(A)\in
\Sigma$ for any $A\in \Sigma$. If $\mu(\varphi^{-1}(A))=0$ for all
$A\in \Sigma$ with $\mu(A)=0$, then $\varphi$ is said to be
nonsingular. This condition means that the measure
$\mu\circ\varphi^{-1}$, defined by
$\mu\circ\varphi^{-1}(A)=\mu(\varphi^{-1}(A))$ for $A\in \Sigma$,
is absolutely continuous with respect to the $\mu$ (it is usually
denoted $\mu\circ\varphi^{-1}\ll \mu$). The Radon-Nikod'ym theorem
ensures the existence of a nonnegative locally integrable function
$h$ on $\Omega$ such that,
$\mu\circ\varphi^{-1}(A)=\int_{A}hd\mu$, $A\in \Sigma$. Any
 nonsingular measurable transformation $\varphi$ induces a linear
operator (composition operator) $C_{\varphi}$ from $L^0(\Omega)$
into itself defined by$$C_{\varphi}(f)(t)=f(\varphi(t)) \ \ \ ;
t\in \Omega, \ \ \ f\in L^0(\Omega),$$ where $L^0(\Omega)$ denotes
the linear space of all equivalence classes of $\Sigma$-measurable
functions on $\Omega$, that is, we identify any two functions that
are equal $\mu$-almost everywhere on $\Omega$. Here the
nonsingularity of $\varphi$ guarantees that the operator
$C_{\varphi}$ is well defined as a mapping from $L^0(\Omega)$ into
itself. If $C_{\varphi}$ maps an Orlicz space $L^{\Phi}(\Omega)$
into itself, then $C_{\varphi}$ is called composition operator on
$L^{\Phi}(\Omega)$. Note that, in this case $C_{\varphi}$ is
bounded.

Let $u:\Omega\rightarrow \mathbb{C}$ be a measurable function on
$\Omega$. Then the rule taking $u$ to $u.f$, is a linear
transformation on $L^0(\Omega)$ and we denote this transformation
by $M_{u}$. In the case that $M_{u}$ is continuous, it is called
multiplication operator induced by $u$.

The composition and multiplication operators received considerable
attention over the past several decades especially on some
measurable function spaces such as $L^P$-spaces, Bergman spaces
and a few ones on Orlicz spaces, such that these operators played
an important role in the study of operators on Hilbert spaces.
\vspace*{0.3cm}

The basic properties of composition and multiplication operators
on measurable function spaces are studied by more mathematicians.
For more details on these operators we refer to Abraham \cite{ab},
 Takagi \cite{ta2}, Axler \cite{alex}, Estaremi and Jabbarzadeh \cite{ej}, Halmos \cite{hal}, Lambert
\cite{lamber}, Singh and Manhas \cite{sin}, Takagi \cite{ta, tam,
taky}, Hudzik and Krbec \cite{h}, Cui, Hudzik, Kumar and
Maligranda \cite{c}, Arora \cite{ar} and some other works. The
multiplication and weighted composition operators are studied on
Orlicz spaces in \cite{ks, skn}. In the case that $\varphi$ is an
N-function, some results on boundedness of composition operators
on Orlicz spaces, are obtained in \cite{ku} (see also \cite{ra}).
As is seen in \cite{sha}, the essential norm plays an interesting
role in the compact problem of concrete operators. Many people
have computed the essential norm of various concrete operators.
For these studies about composition operators, we refer to
\cite{ro, tam, zh}. The question of actually calculating the norm
and essential norm of a composition and multiplication operators
on Orlicz spaces is not a trivial one. In spite of the
difficulties associated with computing the essential norm exactly,
it is often possible to find upper and lower bound for the
essential norm under certain conditions.

In this paper, we are going to present some assertions about
boundedness, compactness and essential norm of multiplication and
composition operators between two Orlicz spaces. In section 2 we
give some necessary and sufficient conditions for boundedness of
composition and multiplication operators between two different
Orlicz spaces. In section 3 we present some necessary and
sufficient conditions for compactness of composition and
multiplication operators between two different Orlicz spaces. Then
in section 4 by using the compactness assertions, that is proved
in section 3, we estimate the essential norm of composition and
multiplication operators.

\section{ \sc\bf Boundedness  }
In this section we present some necessary and sufficient conditions
for boundedness of multiplication and composition operators from
$L^{\Phi_1}(\Omega)$ into $L^{\Phi_2}(\Omega)$.

\vspace*{0.3cm} {\bf Theorem 2.1.} Let $(\Omega, \Sigma, \mu)$ be
a $\sigma$-finite nonatomic measure space and
$\varphi:\Omega\rightarrow \Omega$ be a surjective nonsingular
measurable transformation. Denote by $h$ the Radon-Nikodym
derivative $\frac{d\mu\circ\varphi^{-1}}{d\mu}$. The the following
conditions are equivalent:\\

(a) The composition operator $C_{\varphi}$ is bounded from
$L^{\Phi_1}(\Omega)$ into $L^{\Phi_2}(\Omega)$.\\

(b) The Orlicz space $L^{\Phi_1}(\Omega)$ is embedded continuously
into the weighted orlicz space $L_{h}^{\Phi_2}(\Omega)$, where
$$L_{h}^{\Phi_2}(\Omega)=\{f:\Omega\rightarrow \mathbb{C}:\exists k>0, I_{\Phi_2,h}(f)\int_{\Omega}h\Phi_2(k|f|)d\mu<\infty\}.$$

(c) There are $a, b>0$ and $g\in L^{1}(\Omega)$  such that
$\Phi_{2}(au)h(t)\leq b\Phi_{1}(u)+g(t)$ for all $u>0$ and $t\in
\Omega\setminus A$ with $\mu(A)=0$.\\

 {\bf Proof.} $a\rightarrow b$. Since $\varphi$ is surjective, then for all  $f\in
 L^{\Phi_1}(\Omega)$ we have
 $$I_{\Phi_2}(C_{\varphi}(f))=\int_{\Omega}\Phi_2(|C_{\varphi}(f)|)d\mu
 =\int_{\varphi(\Omega)}h\Phi_2(|f|)d\mu$$$$=\int_{\Omega}h\Phi_2(|f|)d\mu=I_{\Phi_2,h}(f).$$
 Suppose that (a) is satisfied, then for every $f\in
 L^{\Phi_1}(\Omega)$ we have
 $$N_{\Phi_2,h}(f)=N_{\Phi_2}(C_{\varphi}(f))\leq \|C_{\varphi}\|N_{\Phi_1}(f).$$
 This implies that the Orlicz space $L^{\Phi_1}(\Omega)$ is embedded continuously
into the weighted orlicz space $L_{h}^{\Phi_2}(\Omega)$. By
[\cite{m}, Th 8.5], it is easy to see that $b\rightarrow c$.

 For $c\rightarrow a$, we suppose that (c) holds, then for every $f\in L^{\Phi_1}(\Omega)$ we have

 $$I_{\Phi_2}(\frac{aC_{\varphi}(f)}{N_{\Phi_1}(f)})=
 \int_{\Omega}\Phi_{2}(\frac{af(t)}{N_{\Phi_1}(f)})h(t)d\mu$$
 $$\leq b\int_{\Omega}\Phi_{1}(\frac{f(t)}{N_{\Phi_1}(f)})d\mu
 +\int_{\Omega}g(t)d\mu\leq b+\int_{\Omega}g(t)d\mu\leq M',$$
 where $M'>1$. This implies that
 $I_{\Phi_2}(\frac{(aC_{\varphi}(f))}{M'N_{\Phi_1}(f)})\leq1$, thus
 $N_{\Phi_2}(C_{\varphi}(f))\leq\frac{M'}{a}N_{\Phi_1}(f)$.\\

\vspace*{0.3cm} {\bf Theorem 2.2.} If
$C_{\varphi}:L^{\Phi_1}(\Omega)\rightarrow L^{\Phi_2}(\Omega)$ is
a linear transformation, then $C_{\varphi}$ is bounded.\\
{\bf Proof.} By applying closed graph theorem, injectivity of
$\Phi_1$ and $\Phi_2$ and the fact that the norm-convergence
implies the $\Phi$-mean-convergence, we conclude that
$C_{\varphi}$ is bounded.\\

\vspace*{0.3cm} {\bf Remark.2.3} By theorem 2.2 we have:
$C_{\varphi}\in B(L^{\Phi_1}, L^{\Phi_2})$ if and only if
$C_{\varphi}(L^{\Phi_1})\subseteq L^{\Phi_2}$. Thus the following
conditions are equivalent:\\

(a) The composition operator $C_{\varphi}$ is bounded from
$L^{\Phi_1}(\Omega)$ into $L^{\Phi_2}(\Omega)$.\\

(b) For every $f\in L^{\Phi_1}(\Omega)$, there exists $\lambda>0$
such that $$\int_{\Omega}h\Phi_{2}(\lambda|f|)d\mu<\infty.$$

 (c) The Orlicz space $L^{\Phi_1}(\Omega)$ is embedded
continuously into the weighted orlicz space
$L_{h}^{\Phi_2}(\Omega)$.\\

(d) There are $a, b>0$ and $g\in L^{1}(\Omega)$  such that
$\Phi_{2}(au)h(t)\leq b\Phi_{1}(u)+g(t)$ for all $u>0$ and $t\in
\Omega\setminus A$ with $\mu(A)=0$.\\

\vspace*{0.3cm} {\bf Theorem 2.4.} Let $u:\Omega\rightarrow
\mathbb{C}$ be a measurable function. Then\\

(a) If there are $M>0$ and $g\in L_{+}^{1}(\Omega)$  such that
$\Phi_{2}(u(x)v)\leq \Phi_{1}(Mv)+g(x)$ for all $v>0$ and $x\in
\Omega\setminus A$ with $\mu(A)=0$, Then
$M_{u}:L^{\Phi_1}(\Omega)\rightarrow L^{\Phi_2}(\Omega)$ is a
bounded operator.\\

(b) If $(\Omega, \Sigma , \mu)$ is non-atomic measure space and
the operator $M_{u}:L^{\Phi_1}(\Omega)\rightarrow
L^{\Phi_2}(\Omega)$ is bounded, then there exist $M>0$ and $g\in
L_{+}^{1}(\Omega)$  such that $\Phi_{2}(u(x)v)\leq
\Phi_{1}(Mv)+g(x)$ for all $v>0$ and $x\in
\Omega\setminus A$ with $\mu(A)=0$.\\

{\bf Proof.} (a) For every $f\in L^{\Phi_1}(\Omega)$  we have
$$I_{\Phi_2}(\frac{uf}{MN_{\Phi_1}(f)})\leq \int_{\Omega}\Phi_{1}(\frac{Mf(t)}{MN_{\Phi_1}(f)})d\mu
 +\int_{\Omega}g(t)d\mu\leq 1+\int_{\Omega}g(t)d\mu\leq M',$$
 where $M'>1$. This implies that
 $I_{\Phi_2}(\frac{(uf)}{MM'N_{\Phi_1}(f)})\leq1$.
Therefore $N_{\Phi_2}(M_u(f))\leq MM' N_{\Phi_1(f)}$, for every
$f\in L^{\Phi_1}(\Omega)$.\\

(b) If the condition is not satisfied, then for every $n\in
\mathbb{N}$ and $g\in L_{+}^{1}(\Omega)$, there exists a
measurable set $F'_n$ of $\Omega$ and some $\alpha_n\in
\mathbb{C}$ such that
$$F'_n=\{x\in \Omega:\Phi_{2}(|u(x)\alpha_n|)> \Phi_1(2^nn^2\alpha_n)+g(x)\}$$
and $F'_n$ is a measurable set of positive measure. Since
$$F'_n\subseteq \{x\in \Omega:\Phi_{2}(|u(x)\alpha_n|)> \Phi_1(2^nn^2\alpha_n)\}=F_n,$$
then $F_n$ is also a measurable set of positive measure. Since
$\mu$ is non-atomic, we can choose a disjoint sequence of
measurable sets $\{E_n\}$ such that $E_n\subseteq F_n$ and
$$\mu(E_n)=\frac{\Phi_2(|\alpha_1|)}{2^n\Phi_1(n^2|\alpha_n|)}.$$
Let $f=\Sigma^{\infty}_{n=1}c_n\chi_{E_n}$, where
$c_n=n|\alpha_n|$. Suppose that$\alpha>0$ and $n_0>\alpha$, then

$$\int_{\Omega}\Phi_1(\alpha f)d\mu=\Sigma^{\infty}_{n=1}\int_{\Omega}\Phi_1(\alpha
c_n)\chi_{E_n}d\mu$$$$ =\Sigma^{n_0}_{n=1}\Phi_1(\alpha
c_n)\mu(E_n)+\Sigma^{\infty}_{n=n_0+1}\frac{\Phi_1(\alpha
c_n)\Phi_2(|\alpha_1|)}{2^n\Phi_1(n^2|\alpha_n|)}$$
$$\leq\Sigma^{n_0}_{n=1}\Phi_1(\alpha
c_n)\mu(E_n)+\Sigma^{\infty}_{n=n_0+1}\frac{\Phi_1(n^2
|\alpha_n|)\Phi_2(|\alpha_1|)}{2^n\Phi_1(n^2|\alpha_n|)}<\infty.$$
If $\alpha>\frac{1}{n_0}$, then

$$\int_{\Omega}\Phi_2(\alpha |u.f|)d\mu\geq \Sigma_{n\geq n_0}\int_{E_n}\Phi_2(\alpha
n|u\alpha_n|)d\mu\geq$$$$ \Sigma_{n\geq
n_0}\int_{E_n}\Phi_2(|u\alpha_n|)d\mu\geq \Sigma_{n\geq
n_0}\Phi_1(2^nn^2|\alpha_n|)\mu(E_n)$$$$\geq \Sigma_{n\geq
n_0}\Phi_2(|\alpha_1)=\infty.$$ This is a contradiction.\\

\vspace*{0.3cm} {\bf Corollary 2.5.} Let $(\Omega, \Sigma, \mu)$
be a non-atomic measure space. Then
$M_{u}:L^{\Phi_1}(\Omega)\rightarrow L^{\Phi_2}(\Omega)$ is a
bounded operator if and only if there exist $M>0$ and $g\in
L_{+}^{1}(\Omega)$ such that $\Phi_{2}(u(x)v)\leq
\Phi_{1}(Mv)+g(x)$ for all $v>0$ and $x\in \Omega\setminus A$ with
$\mu(A)=0$.\\

\vspace*{0.3cm} {\bf Theorem 2.6.} If
$M_{u}:L^{\Phi_1}(\Omega)\rightarrow L^{\Phi_2}(\Omega)$ is a
linear transformation, then $M_{u}$ is bounded.\\
{\bf Proof.} By applying closed graph theorem, injectivity of
$\Phi_1$ and $\Phi_2$ and the fact that the norm-convergence
implies the $\Phi$-mean-convergence, we conclude that
$C_{\varphi}$ is bounded.

\section{ \sc\bf Compactness }
In this section we present some necessary and sufficient condition
for composition and multiplication operators to be compact.
 Recall that an atom of the measure $\mu$ is
an element $A\in\Sigma$ with $\mu(A)>0$ such that for each
$F\in\Sigma$, if $F\subseteq A$ then either $\mu(F)=0$ or
$\mu(F)=\mu(A)$. A measure space $(\Omega,\Sigma,\mu)$ with no
atoms is called non-atomic measure space. It is well-known fact
that every $\sigma$-finite measure space $(\Omega, \Sigma,\mu)$
can be partitioned uniquely as $\Omega. =B\cup\{A_j:
j\in\mathbb{N}\}$, where $\{A_j\}_{j\in\mathbb{N}}$ is a countable
collection of pairwise disjoint atoms and $B\in\Sigma$, being
disjoint from each $A_j$, is non-atomic (see \cite{z}). Since
$\Sigma$ is $\sigma$-finite, so $a_j:=\mu(A_j)<\infty$, for all
$j\in\mathbb{N}$. A bounded linear operator $T:E\rightarrow E$
(where $E$ is a Banach space) is called compact, if $T(B_1)$ has
compact closure, where $B_1$ denotes the closed unit ball of $E$.\\

 \vspace*{0.3cm} {\bf Theorem 3.1.} Let $T=C_{\varphi}$ be bounded from $L^{\Phi_{1}}(\Omega, \Sigma,
\mu)$ to $ L^{\Phi_{2}}(\Omega, \Sigma, \mu)$. Then $C_{\varphi}$
is compact if and only if $N_{\varepsilon}=\{x\in
\Omega:\Phi_{2}(|\alpha|)h(t)>\Phi_{1}(\varepsilon |\alpha|),
\alpha\in \mathbb{C}\}$ consists of finitely many atoms, for all
$\varepsilon>0$.\\

{\bf Proof.} Let $\varepsilon>0$ and
$N_{\varepsilon}=\cup^{n}_{i=1}C_{n}$ consists of finitely many
atoms. Put
$T_{\varepsilon}=C_{\varphi}M_{\chi_{N_{\varepsilon}}}$. Since
$\Sigma-$measurable functions are constant on $\Sigma-$atoms and
$(\Omega, \Sigma, \mu)$ is $\sigma-$finite, then
$M_{\chi_{N_{\varepsilon}}}$ is a compact operator on
$L^{\Phi_{1}}(\Omega, \Sigma, \mu)$.
 Thus the operator $T_{\varepsilon}=C_{\varphi}M_{\chi_{N_{\varepsilon}}}$ is compact from $L^{\Phi_{1}}(\Omega, \Sigma,
\mu)$ to $ L^{\Phi_{2}}(\Omega, \Sigma, \mu)$.
 Hence for every $f\in
L^{\Phi_{1}}(\Omega, \Sigma, \mu)$

$$\int_{\Omega}\Phi_{2}(\frac{(T-T_{\varepsilon})(f)}{\varepsilon N_{\Phi_{1}}(f)} )d\mu
=\int_{\Omega}\Phi_{2}(\frac{C_{\varphi}(\chi_{\Omega\setminus N_{\varepsilon}}f)}{\varepsilon N_{\Phi_{1}}(f)})
d\mu$$
$$=\int_{\Omega\setminus N_{\varepsilon}}h\Phi_{2}(\frac{f}{\varepsilon N_{\Phi_{1}}(f)})
d\mu\leq\int_{\Omega\setminus N_{\varepsilon}}\Phi_{1}(\frac{\varepsilon f}{\varepsilon N_{\Phi_{1}}(f)}) d\mu$$
$$=\int_{\Omega\setminus N_{\varepsilon}}\Phi_{1}(\frac{f}{N_{\Phi_{1}}(f)}) d\mu\leq1.$$
This implies that $N_{\Phi_{2}}(Tf-T_{\varepsilon}f)\leq
\varepsilon N_{\Phi_{1}}(f)$. Thus $T$ is compact.\\

Conversely, suppose there exists $\varepsilon>0$ such that $N_{\varepsilon}$ consists of infinitely many atoms or a non-atomic subset of positive measure. In both cases we can find a sequence $\{B_{n}\}_{n\in \mathbb{N}}\}$ of disjoint measurable subsets of $N_{\varepsilon}$ with $0<\mu(B_{n})<\infty$. Put $f_{n}=\frac{\chi_{B_{n}}}{N_{\Phi_{1}}(\chi_{B_{n}})}$. Hence

$$\int_{\Omega}\Phi_{1}(\frac{\varepsilon f_{n}}{N_{\Phi_{2}}(f_{n}\circ\varphi)})d\mu\leq \int_{\Omega}h\Phi_{2}(\frac{f_{n}}{N_{\Phi_{2}}(f_{n})})d\mu$$$$=\int_{\Omega}\Phi_{2}(\frac{f_{n}\circ\varphi}{N_{\Phi_{2}}(f_{n}\circ\varphi)})d\mu\leq1.$$

So $\varepsilon=N_{\Phi_{1}}(\varepsilon f_{n})\leq N_{\Phi_{2}}(f_{n}\circ\varphi)$. Since $B_{n}$'s are disjoint. for $n\neq m$
$$N_{\Phi_{2}}(f_{n}\circ\varphi-f_{m}\circ\varphi)\geq N_{\Phi_{2}}(f_{n}\circ\varphi)\geq \varepsilon.$$
So $\{f_{n}\circ\varphi\}_{n\in \mathbb{N}}$ has no convergent
subsequence. This mean's $T=C_{\varphi}$ can not be compact.\\

\vspace*{0.3cm} {\bf Corollary 3.2.} If $(\Omega, \sigma, \mu)$ is
nonatomic measure space, then there is not nonzero compact
composition operator between  $L^{\Phi_{1}}(\Omega, \Sigma, \mu)$
and $ L^{\Phi_{2}}(\Omega, \Sigma, \mu)$.\\

 \vspace*{0.3cm} {\bf Theorem 3.3.} Let $M_{u}$ be bounded from $L^{\Phi_{1}}(\Omega, \Sigma,
\mu)$ to $ L^{\Phi_{2}}(\Omega, \Sigma, \mu)$. Then $M_{u}$ is
compact if and only if $N_{\varepsilon}(u)=\{x\in
\Omega:\Phi_{2}(|u(x)\alpha|)>\Phi_{1}(\varepsilon |\alpha|),
\alpha\in \mathbb{C}\}$ consists of finitely many atoms, for all
$\varepsilon>0$.\\

{\bf Proof.} Let $\varepsilon>0$ and
$N_{\varepsilon}=N_{\varepsilon}(u)=\cup^{n}_{i=1}C_{n}$ consists
of finitely many atoms. Put
$u_{\varepsilon}=u\chi_{N_{\varepsilon}}$ and
$T_{\varepsilon}=M_{u_{\varepsilon}}$. Since $\Sigma-$measurable
functions are constant on $\Sigma-$atoms and $(\Omega, \Sigma,
\mu)$ is $\sigma-$finite. We have

$$M_{u_{\varepsilon}}(f)=\sum^{n}_{i=1}u(C_{i})f(C_{i})\chi_{C_{i}}\in \{\sum^{n}_{i=1}\alpha_{i}\chi_{C_{i}}:\alpha_{i}\in \mathbb{C}\}\subseteq L^{\Phi_{2}}(\Omega, \Sigma, \mu).$$

Thus $M_{u_{\varepsilon}}$ is finite rank. Hence for every $f\in
L^{\Phi_{1}}(\Omega, \Sigma, \mu)$

$$\int_{\Omega}\Phi_{2}(\frac{(u-u_{\varepsilon})f}{\varepsilon N_{\Phi_{1}}(f)} )d\mu=\int_{\Omega\setminus N_{\varepsilon}}\Phi_{2}(\frac{uf}{\varepsilon N_{\Phi_{1}}(f)}) d\mu$$
$$\leq\int_{\Omega\setminus N_{\varepsilon}}\Phi_{1}(\frac{\varepsilon f}{\varepsilon N_{\Phi_{1}}(f)}) d\mu=\int_{\Omega\setminus N_{\varepsilon}}\Phi_{1}(\frac{f}{N_{\Phi_{1}}(f)}) d\mu\leq1.$$
This implies that $N_{\Phi_{2}}(M_{u}f-M_{u_{\varepsilon}}f)\leq
\varepsilon N_{\Phi_{1}}(f)$. Thus $M_{u}$ is compact.\\

Conversely, suppose there exists $\varepsilon>0$ such that
$N_{\varepsilon}$ consists of infinitely many atoms or a
non-atomic subset of positive measure. In both cases we can find a
sequence $\{B_{n}\}_{n\in \mathbb{N}}\}$ of disjoint measurable
subsets of $N_{\varepsilon}$ with $0<\mu(B_{n})<\infty$. Put
$f_{n}=\frac{\chi_{B_{n}}}{N_{\Phi_{1}}(\chi_{B_{n}})}$. Hence

$$\int_{\Omega}\Phi_{1}(\frac{\varepsilon f_{n}}{N_{\Phi_{2}}(uf_{n})})d\mu\leq \int_{\Omega}\Phi_{2}(\frac{uf_{n}}{N_{\Phi_{2}}(uf_{n})})d\mu\leq1.$$

So $\varepsilon=N_{\Phi_{1}}(\varepsilon f_{n})\leq N_{\Phi_{2}}(u
f_{n})$. Since $B_{n}$'s are disjoint. for $n\neq m$
$N_{\Phi_{2}}(u f_{n}-uf_{m})\geq N_{\Phi_{2}}(u f_{n})\geq
\varepsilon$. So $\{uf_{n}\}_{n\in \mathbb{N}}$ has no convergent
subsequence. This mean's $M_{u}$ can not be compact.\\

\vspace*{0.3cm} {\bf Corollary 3.4.} If $(\Omega, \sigma, \mu)$ is
nonatomic measure space, then there is not nonzero compact
multiplication operator between  $L^{\Phi_{1}}(\Omega, \Sigma,
\mu)$ and $ L^{\Phi_{2}}(\Omega, \Sigma, \mu)$.

\section{ \sc\bf Essential norm }
Let $\mathfrak{B}$ be a Banach space and $\mathcal{K}$ be the set
of all compact operators on $\mathfrak{B}$. For $T\in
L(\mathfrak{B})$, the Banach algebra of all bounded linear
operators on $\mathfrak{B}$ into itself, the essential norm of $T$
means the distance from $T$ to $\mathcal{K}$ in the operator norm,
namely $\|T\|_e =\inf\{\|T - S\| : S \in\mathcal{K}\}$. Clearly,
$T$ is compact if and only if $\|T\|_e= 0$. As is seen in
\cite{sha}, the essential norm plays an interesting role in the
compact problem of concrete operators.

\vspace*{0.3cm} {\bf Theorem 4.1.} Let
$\varphi:\Omega\rightarrow\Omega$
 be nonsingular measurable transformation
 and let $T=C_{\varphi}: L^{\Phi_{1}}(\Omega, \Sigma, \mu)\rightarrow L^{\Phi_{2}}(\Omega, \Sigma, \mu)$.
 If
$\beta_1=\inf\{\varepsilon>0:N_{\varepsilon}$ consists of finitely
many atoms$\}$. Then\\

(a) $\|C_{\varphi}\|_{e}\leq \beta_1$.\\

(b) Let $\Phi_{1}\in \bigtriangleup_{2}$ and
$\mu(C_{n})\rightarrow 0$ or $\{\mu(C_{n})\}_{n\in \mathbb{N}}$
has no convergent subsequence. Then $\beta_1\leq
\|C_{\varphi}\|_{e}$.\\

{\bf Proof.} (a) Let $\varepsilon>0$. Then
$N_{\varepsilon+\beta_1}$ consist of finitely many atoms.
 Put $T_{\varepsilon+\beta_1}=C_{\varphi}M_{\chi_{N_{\varepsilon+\beta_1}}}$.
 So $T_{\varepsilon+\beta_1}$ is compact. Also for $f\in L^{\Phi_{1}}(\Omega, \Sigma, \mu)$

$$\int_{\Omega}\Phi_{2}(\frac{(T-T_{\varepsilon+\beta_1})(f)}{(\varepsilon+\beta_1) N_{\Phi_{1}}(f)} )d\mu
=\int_{\Omega}\Phi_{2}(\frac{C_{\varphi}(\chi_{\Omega\setminus
N_{\varepsilon+\beta_1}}f)}{(\varepsilon+\beta_1)
N_{\Phi_{1}}(f)}) d\mu$$
$$=\int_{\Omega\setminus N_{\varepsilon+\beta_1}}h\Phi_{2}(\frac{f}{(\varepsilon+\beta_1) N_{\Phi_{1}}(f)})
d\mu\leq\int_{\Omega\setminus
N_{\varepsilon+\beta_1}}\Phi_{1}(\frac{(\varepsilon+\beta_1)
f}{(\varepsilon+\beta_1) N_{\Phi_{1}}(f)})
d\mu$$$$=\int_{\Omega\setminus
N_{\varepsilon+\beta_1}}\Phi_{1}(\frac{f}{N_{\Phi_{1}}(f)})
d\mu\leq1.$$ This implies that
$N_{\Phi_{2}}(Tf-T_{\varepsilon+\beta_1}f)\leq
(\varepsilon+\beta_1) N_{\Phi_{1}}(f)$. Hence
$$\|T\|_{e}\leq \|T-T_{\varepsilon+\beta_1}\|\leq \varepsilon+\beta_1.$$
Thus $\|T\|_{e}\leq \beta_1$.\\

(b) Let $0<\varepsilon<\beta_1$. Then by definition,
$N_{\beta_1-\varepsilon}(u)$ contains infinitely many atoms or a
non- atomic subset of positive measure. If
$N_{\beta_1-\varepsilon}(u)$ consists a non- atomic subset, then
we can find a sequence $\{B_{n}\}_{n\in \mathbb{N}}$ such that
$\mu(B_{n})<\infty$ and $\mu(B_{n})\rightarrow 0$. Put
$f_{n}=\frac{\chi_{B_{n}}}{N_{\Phi_{1}}(\chi_{B_{n}})}$, then for
every $A\in \Sigma$ with $0<\mu(A)<\infty$ we have

$$\int_{\Omega}f_{n}\chi_{A}d\mu=\mu(A\cap B_{n})\Phi^{-1}_{1}(\frac{1}{\mu(B_{n})})\leq \frac{\Phi^{-1}_{1}(\frac{1}{\mu(B_{n})})}{\frac{1}{\mu(B_{n})}}
\rightarrow0.$$ when $n\rightarrow \infty$. Also, if
$N_{\beta_1-\varepsilon}(u)$ consists infinitely many atoms
$\{C'_{n}\}_{n\in \mathbb{N}}$. We set
$f_{n}=\frac{\chi_{C'_{n}}}{N_{\Phi_{1}}(\chi_{C'_{n}})}$. Then
for every $A\in \Sigma$ with $0<\mu(A)<\infty$ we have
$$\int_{\Omega}f_{n}\chi_{A}d\mu=\mu(A\cap C'_{n})\Phi^{-1}_{1}(\frac{1}{\mu(C'_{n})}).$$
If $\{\mu(C_{n})\}_{n\in \mathbb{N}}$ has no convergent
subsequence, then there exists $n_{0}$ such that for $n>n_{0}$,
$\mu(A\cap C_{n}')=0$ and if $\mu(C_{n})\rightarrow 0$ then
$\mu(C_{n}')\rightarrow 0$. Thus
$\int_{\Omega}f_{n}\chi_{A}d\mu=\mu(A\cap
C'_{n})\Phi^{-1}_{1}(\frac{1}{\mu(C'_{n})})\rightarrow0$ in both
cases. These imply that $f_{n}\rightarrow0$ weakly. So
$$\int_{\Omega}\Phi_{1}(\frac{(\beta_1-\varepsilon)f_{n}}{N_{\Phi_{2}}(f_{n}\circ\varphi)})d\mu
\leq\int_{\Omega}\Phi_{2}(\frac{f_{n}\circ\varphi}{N_{\Phi_{2}}(f_{n}\circ\varphi)})d\mu.$$

Thus $N_{\Phi_{2}}(C_{\varphi}(f_{n}\circ\varphi))\geq
\beta_1-\varepsilon$.\\

Also, there exists compact operator $T\in L(L^{\Phi_{1}}(\Omega,
\Sigma, \mu), L^{\Phi_{1}}(\Omega, \Sigma, \mu))$ such that
$\|C_{\varphi}\|_{e}\geq\|C_{\varphi}-T\|-\varepsilon$. Hence
$N_{\Phi_{2}}(Tf_{n})\rightarrow o$ and so there exists $N>0$ such
that for each $n>N$, $N_{\Phi_{2}}(Tf_{n})\leq \varepsilon$. So

$$\|C_{\varphi}\|_{e}\geq\|C_{\varphi}-T\|-\varepsilon\geq|N_{\Phi_{2}}(f_{n}\circ\varphi)-N_{\Phi_{2}}(Tf_{n})|\geq\beta_1-\varepsilon-\varepsilon,$$
thus we conclude that $\|C_{\varphi}\|_{e}\geq\beta_1$.\\

\vspace*{0.3cm} {\bf Theorem 4.2.} Let
$u:\Omega\rightarrow\mathbb{C}$ be $\Sigma-$measurable and Let
$M_{u}: L^{\Phi_{1}}(\Omega, \Sigma, \mu)\rightarrow
L^{\Phi_{2}}(\Omega, \Sigma, \mu)$. If
$\beta_2=\inf\{\varepsilon>0:N_{\varepsilon}$ consists of finitely
many atoms$\}$. Then\\

(a) $\|M_{u}\|_{e}\leq \beta_2$.\\

(b) Let $\Phi_{1}\in \bigtriangleup_{2}$ and
$\mu(C_{n})\rightarrow 0$ or $\{\mu(C_{n})\}_{n\in \mathbb{N}}$
has no convergent subsequence. Then $\beta_2\leq \|M_{u}\|_{e}$.\\

{\bf Proof.} (a) Let $\varepsilon>0$. Then
$N_{\varepsilon+\beta_2}$ consist of finitely many atoms. Put
$u_{\varepsilon+\beta_2}=u\chi_{N_{\varepsilon+\beta_2}}$ and
$M_{u_{\varepsilon+\beta_2}}$. So $M_{u_{\varepsilon+\beta_2}}$ is
finite rank and so compact. Also for $f\in L^{\Phi_{1}}(\Omega,
\Sigma, \mu)$

$$\int_{\Omega}\Phi_{2}(\frac{(u-u_{\varepsilon+\beta_2})f}{(\varepsilon+\beta_2) N_{\Phi_{1}}(f)} )d\mu
=\int_{\Omega\setminus
N_{\varepsilon+\beta_2}}\Phi_{2}(\frac{uf}{(\varepsilon+\beta_2)
N_{\Phi_{1}}(f)}) d\mu$$
$$\leq\int_{\Omega\setminus N_{\varepsilon
+\beta_2}}\Phi_{1}(\frac{f}{N_{\Phi_{1}}(f)}) d\mu\leq1.$$

Hence $N_{\Phi_{2}}(uf-u_{\varepsilon+\beta_2}f)\leq
(\varepsilon+\beta_2)N_{\Phi_{1}}(f)$ and so

$$\|M_{u}\|_{e}\leq \|M_{u}-M_{u_{\varepsilon+\beta_2}}\|\leq \varepsilon+\beta_2.$$
Thus $\|M_{u}\|_{e}\leq \beta_2$.\\

(b) Let $0<\varepsilon<\beta_2$. Then by definition,
$N_{\beta_2-\varepsilon}(u)$ contains infinitely many atoms or a
non- atomic subset of positive measure. If
$N_{\beta_2-\varepsilon}(u)$ consists a non- atomic subset, then
we can find a sequence $\{B_{n}\}_{n\in \mathbb{N}}$ such that
$\mu(B_{n})<\infty$ and $\mu(B_{n})\rightarrow 0$. Put
$f_{n}=\frac{\chi_{B_{n}}}{N_{\Phi_{1}}(\chi_{B_{n}})}$, then for
every $A\in \Sigma$ with $0<\mu(A)<\infty$ we have

$$\int_{\Omega}f_{n}\chi_{A}d\mu=\mu(A\cap B_{n})\Phi^{-1}_1(\frac{1}{\mu(B_{n})})\leq \frac{\Phi^{-1}_1(\frac{1}{\mu(B_{n})})}{\frac{1}{\mu(B_{n})}}
\rightarrow0.$$ when $n\rightarrow \infty$. Also, if
$N_{\beta_2-\varepsilon}(u)$ consists infinitely many atoms
$\{C'_{n}\}_{n\in \mathbb{N}}$. We set
$f_{n}=\frac{\chi_{C'_{n}}}{N_{\Phi_{1}}(\chi_{C'_{n}})}$. Then
for every $A\in \Sigma$ with $0<\mu(A)<\infty$ we have
$$\int_{\Omega}f_{n}\chi_{A}d\mu=\mu(A\cap C'_{n})\Phi^{-1}_1(\frac{1}{\mu(C'_{n})}).$$
If $\{\mu(C_{n})\}_{n\in \mathbb{N}}$ has no convergent
subsequence, then there exists $n_{0}$ such that for $n>n_{0}$,
$\mu(A\cap C_{n}')=0$ and if $\mu(C_{n})\rightarrow 0$ then
$\mu(C_{n}')\rightarrow 0$. Thus
$\int_{\Omega}f_{n}\chi_{A}d\mu=\mu(A\cap
C'_{n})\Phi^{-1}_1(\frac{1}{\mu(C'_{n})})\rightarrow0$ in both
cases. These imply that $f_{n}\rightarrow0$ weakly. So
$$\int_{\Omega}\Phi_{1}(\frac{(\beta-\varepsilon)f_{n}}{N_{\Phi_{2}}(uf_{n})})d\mu
\leq\int_{\Omega}\Phi_{2}(\frac{uf_{n}}{N_{\Phi_{2}}(uf_{n})})d\mu.$$

Thus $N_{\Phi_{2}}(uf_{n})\geq \beta_2-\varepsilon$.

Also, there exists compact operator $T\in L(L^{\Phi_{1}}(\Omega,
\Sigma, \mu), L^{\Phi_{2}}(\Omega, \Sigma, \mu))$ such that
$\|M_{u}\|_{e}\geq\|T-M_{u}\|-\varepsilon$. Hence
$N_{\Phi_{2}}(Tf_{n})\rightarrow o$ and so there exists $N>0$ such
that for each $n>N$, $N_{\Phi_{2}}(Tf_{n})\leq \varepsilon$. So

$$\|M_{u}\|_{e}\geq\|M_{u}-T\|-\varepsilon\geq|N_{\Phi_{2}}(uf_{n})-N_{\Phi_{2}}(Tf_{n})|\geq\beta_2-\varepsilon-\varepsilon,$$
thus we conclude that $\|M_{u}\|_{e}\geq\beta_2$.\\

\vspace*{0.3cm} {\bf Corollary 4.3.} If $\mu(\Omega)<\infty$ and
$\Phi_{1}\in \bigtriangleup_{2}$.
 Then\\

 (a) $\|C_{\varphi}\|_{e}=\beta_1$.\\

 (b) $\|M_{u}\|_{e}=\beta_2$.\\

 \vspace*{0.3cm} {\bf Corollary 4.4.} If $\mu(\Omega)<\infty$ and
$\Phi_{1}\in \bigtriangleup_{2}$.
 Then\\

 (a) $C_{\varphi}$ is compact if and only if $\beta_1=0$.\\

 (b) $M_u$ is compact if and only if $\beta_2=0$.\\

\vspace*{0.3cm} {\bf Example 4.5.}  (a) Let $\Omega=\mathbb{N}$,
$\mu$ be counting measure on $\Omega$ and $\varphi$ be injective
transformation on $\Omega$. If we set $\Phi_1(n)=\frac{n^3}{3}$
and $\Phi_2(n)=\frac{n^2}{2}$, for all $n\in \mathbb{N}$. Then by
using theorem 3.1 the operator $C_{\varphi}$ is compact from
$L^{\Phi_{1}}(\mathbb{N}, \Sigma, \mu)$ into
$L^{\Phi_{2}}(\mathbb{N}, \Sigma, \mu)$. Also, $C_{\varphi}$ is
not compact from $L^{\Phi_{2}}(\mathbb{N}, \Sigma, \mu)$ into
$L^{\Phi_{1}}(\mathbb{N}, \Sigma, \mu)$. Also, if
$u(n)=\frac{n^2}{n+1}$, then by theorem 3.3 the multiplication
operator $M_u$ is not compact from $L^{\Phi_{1}}(\mathbb{N},
\Sigma, \mu)$ into $L^{\Phi_{2}}(\mathbb{N}, \Sigma, \mu)$ and if
$u(n)=\frac{1}{n^2}$, then by theorem 3.3 the multiplication
operator $M_u$ is compact from $L^{\Phi_{1}}(\mathbb{N}, \Sigma,
\mu)$ into $L^{\Phi_{2}}(\mathbb{N}, \Sigma, \mu)$.\\

(b) Let $\Omega=[0, 1)\cup (\mathbb{N}-\{1\})$, where $\mathbb{N}$
is the set of natural numbers. Let $\mu$ be the Lebesque measure
on $[0, 1)$ and $\mu(\{n\})=1$, if $n\in(\mathbb{N}-\{1\})$. If we
set $\Phi_1(x)=e^x-x-1$, $\Phi_2(n)=\frac{x^5}{5}$  and
$u(x)=x^2+2$ for $x\in \Omega$, then then by theorem 3.3 the
multiplication operator $M_u$ is not compact from
$L^{\Phi_{1}}(\mathbb{N}, \Sigma, \mu)$ into
$L^{\Phi_{2}}(\mathbb{N}, \Sigma, \mu)$.

\vspace*{0.3cm} {\bf ACKNOWLEDGEMENT} The author would like to
thank the professor Romesh Kumar  for very helpful comments and
valuable suggestions.


\begin{thebibliography}{99}

\bibitem{ab}  M.B. Abrahamse , Multiplication Operators, Lecture Notes in Mathematics, {\bf693}, (1978), 17-36.

\bibitem{alex}  A. Axler, Multiplication Operators on Bergman Spaces, J.
Reine,. Angew Math. {\bf336} (1982), 26-44.

\bibitem{ar}
S. C. Arora, Gopal Datt and Satish Verma, Composition operators on
Lorentz spaces, Bull. Austral. Math. Soc. {\bf 76} (2007),
205-214.

\bibitem{cyh} Cui, Yunan; Hudzik, Henryk; Kumar, Romesh; Maligranda, Lech
Composition operators in Orlicz spaces. J. Aust. Math. Soc. 76
(2004), no. 2, 189�206.

\bibitem{doug} R.G. Douglas, Banach Algebra Techniques in Operator Theory, Academic Press,
New York,  1972.

\bibitem{ej} Y. Estaremi and M.R. Jabbarzadeh, Weighted lambert type operators on
$L^{p}$-spaces, Oper. and Matrices. {\bf 7} (2013), 101-116.

\bibitem{hal} P.R. Halmos, A Hilbert
Space Problem Book, Van  Nostr and, Princeton, N.J., 1961.

\bibitem{h}
H. Hudzik and M. Krbec,  On non-effective weights in Orlicz
spaces. Indag. Math. (N.S.) {\bf 18} (2007),  215-231.

\bibitem{c}
Y. Cui, H. Hudzik, R. Kumar and L. Maligranda,  Composition
operators in Orlicz spaces, J. Aust. Math. Soc. {\bf 76} (2004),
189-206.

\bibitem{ks} B.S.  Komal  AND  S. Gupta, Multiplication operators between Orlicz
spaces, Integral equation and operator theory. {\bf41} (2001),
324-330.

\bibitem{ku} R. Kumar, 'Composition operators on Orlicz spaces', Integral Equations Operator Theory
{\bf29} (1997), 17–22.


\bibitem{m} J. Musielak, Orlicz spaces and modular spaces, Lecture Notes in
Math.1034 (Springer, Berlin, 1983).

%\bibitem{lambe}
%A. Lambert, Localising sets for sigma-algebras and related point
%transformations, Proc. Roy. Soc. Edinburgh Sect. A {\bf 118}
%(1991), 111-118.

\bibitem{lamber}
A. Lambert, Operator algebras related to measure preserving
transformations of finite order, Rocky Mountain J. Math. {\bf 14}
(1984), 341-349

\bibitem{ra} M. M. Rao, 'Convolutions of vector fields–II: random walk models', Nonlinear Anal, Theory
Methods Appl. {\bf47} (2001), 3599-3615.

\bibitem{ro} Montes-Rodr´ýguez, A. The essential norm of a composition operator on Bloch spaces, Pacific J. Math. 188, 339-351
(1999).

\bibitem{sin} R. K. Singh and J. S. Manhas, Composition operators on function
spaces, North Holland Math. Studies 179, Amsterdam 1993.

\bibitem{skn} S. Gupta, B. S Komal and N. Suri, Weighted
composition operators on Orlicz spaces, Int. J. Contemp. Math.
Sciences. {\bf 1}, 11-20 (2010).

\bibitem{sha}
J. H. Shapiro, The essential norm of a composition operator,
Analls of Math. {\bf 125} (1987), 375-404.

\bibitem{ta}
H. Takagi, Compact weighted composition operators on $L^p$, Proc.
Amer. Math. Soc. {\bf 116} (1992), 505-511.

\bibitem{ta2}
 H.  Takagi,  FredhoIm Weighted Composition  Operators, Integral Equations and Operator Theory {\bf 16} (1993),
267-276.

\bibitem{tam} Takagi, T., Miura, T. and Takahasi, Sin-Ei: Essential norm and stability constants of weighted composition operators
on $C(X)$, Bull. Korean Math. Soc. {\bf 40}, 583-591 (2003).

\bibitem{taky}
H. Takagi and K. Yokouchi, Multiplication and composition operators
between two $L^p$-spaces, Contemporary Math. {\bf 232}(1999),
321-338.

\bibitem{z}
A. C. Zaanen, Integration, 2nd ed., North-Holland, Amsterdam,
1967.

\bibitem{zh} Zheng, L.: The essential norms and spectra of composition operators on $H^{\infty}$, Pacific J. Math. {\bf 203}, 503-510 (2002).

\end{thebibliography}
\end{document}